\documentclass{amsart}

\usepackage{amsmath, amsthm}
\usepackage{amssymb, latexsym}
\usepackage{amsfonts}

 \usepackage{graphicx}
\usepackage[dvipsnames]{xcolor}

\usepackage{isomath}

\usepackage{tkz-graph}

\usepackage{calrsfs}

\newtheorem{theorem}[equation]{Theorem}

\newtheorem{Lemma}[equation]{Lemma}

\newtheorem{rema}[equation]{Remark}
\newtheorem{Remark}[equation]{Remark}

\newtheorem{coro}[equation]{Corollary}
\newtheorem{prop}[equation]{Proposition}

\numberwithin{equation}{section}
   
  \DeclareMathSymbol{\upGamma}{\mathalpha}{operators}{0}

\begin{document}

\title[Automorphism groups of designs with $\lambda=1$]
{Automorphism groups of designs with $\lambda=1$
 }

       \author{William  M. Kantor
      }
       \address{U. of Oregon,
       Eugene, OR 97403
        \      and  \
       Northeastern U., Boston, MA 02115}
  \email{kantor@uoregon.edu}

\begin{abstract}
\vspace{-4pt}
If $G$ is a finite group and $k =q>2$ or $k=q+1$ for   a prime power $q$
 then, for infinitely many integers $v$,  there is a 
$2$--$(v,k,1)$-design ${\bf D}$ for which  ${\rm Aut} {\bf D}\cong G$.

 \end{abstract}

\maketitle
 
\section{Introduction
   }
 \label{Introduction}
  
Starting with Frucht's theorem on graphs \cite{Frucht}, there have been many papers  proving that any finite group is isomorphic to the full automorphism group of 
some specific type of combinatorial object.
  Babai surveyed this topic \cite{Ba1}, and in  \cite[p.~8]{Ba1} stated that 
  in \cite{Ba2}
  he had  proved  that 
  2-designs with $\lambda=1$ are such objects when 
  $k =q>2$ or $k=q+1$ for a prime power~$q$.
  (The case of Steiner triple systems was handled in \cite{Me}.)
  The purpose of this note is to provide a proof of Babai's result\footnote{This theorem was proved before I knew of Babai's result.}:

\begin{theorem}
\label{design construction3}
Let $G$ be a finite group and $q$ a prime power.
 \begin{itemize}
\item  [\rm(i)]There are   infinitely many integers $v $ such that there is a
$2$--$(v,q+1,1)$-design ${\bf D}$ for which  ${\rm Aut} {\bf D}\cong G$.

\item [\rm(ii)]
If  $q>2$ then there are   infinitely many integers $v $ such that there is a 
\mbox{$2$--$(v,q,1)$-design} ${\bf D}$ for which  ${\rm Aut} {\bf D}\cong G$.
\end{itemize}
\end{theorem}

Parts of our proof mimic    \cite[Sec.~5]{DK} and \cite[Sec.~4]{Ka}, but the present situation   is much simpler.
We modify a   small number of subspaces of a projective 
or affine space in such a way that the projective 
or affine space can be recovered   from the resulting design by elementary
geometric arguments.  Further geometric 
arguments determine the automorphism group. 

%\smallskip\smallskip
Section   \ref{Concluding remarks}  contains further properties of the design {\bf D} in the theorem,
 some of  which are needed  in  future research \cite{DoK}.

   \medskip
\noindent
\emph{Notation}:
 We  use standard permutation group notation, such as 
 $x^\pi$ for the image of a point $x$ under a permutation $\pi$ and
 $g^h=h ^{-1}gh $ for conjugation.
 The group of automorphisms of a projective space $Y={\rm PG}(V)$ defined by a vector space $V$ is denoted by ${\rm P\upGamma L}(V)={\rm P\upGamma L}(Y)$; this is induced by the group ${\rm \upGamma L}(V)$   of
 invertible semilinear transformations on $V$.  Also ${\rm A\upGamma L}(V)$ denotes the 
 group of automorphisms of  the affine  space  ${\rm AG}(V)$ defined by   $V$.
 
\section{A simple projective  construction}
\label{Construction}
Let $G$ be a finite group.  For suitable $n$  let 
$\Gamma$ be  a simple, undirected, connected  graph on $\{1,\ldots ,n\}$  such that
${\rm Aut} \Gamma  \cong G$ and $G$ acts semiregularly on the vertices (as in \cite{Frucht}).
All sufficiently large $n$ satisfy  these conditions; we always assume that $n\ge6$.

Let $K={\bf F}_q \subset F={\bf F}_{q^4}$,  and let $\theta$ generate $F^*$.  
Let $V_F$ be an $n$-dimensional  vector space over   $F$,
with basis  $ v_1,\ldots ,v_n $.
View $G$ as acting on $V_F$, permuting   
$\{ v_1,\ldots ,v_n \}$ as it does $\{1,\ldots ,n\}$.
View $V_F$  as a vector space $V$ over $K$.   
  If $Y$ is a set of points of ${\bf  P}={\rm PG}(V)$
   then $\langle Y\rangle$ denotes the smallest subspace
   of {\bf P}   containing $Y$.
   
   We will modify the 
point-line design ${\rm PG}_1(V)$ of  ${\bf  P} $, using  nonisomorphic designs 
$  \Delta _1$ and  $\Delta _2 $ 
whose parameters are those of   ${\rm PG}_1(K^4)={\rm PG}_1(3,q)$   but are not isomorphic to 
that design,  chosen so that   ${\rm Aut} \Delta _1$
fixes a point (Proposition~\ref{needed designs}).  

Our design $\bf D$ has 
the set  $\mathfrak P$ of points of  {\bf  P}
as its set of points.  Most   blocks of $\bf D$  
are lines of  {\bf  P}, with the following exceptions involving  some of 
the subspaces $Fv, $ $0\ne v\in V$, viewed as subsets of $\mathfrak P$.  
For orbit representatives $i$ and $ij$ 
of $G$ on the vertices and edges  of $\Gamma$,  
\begin{enumerate}
\item [(I)]   replace the set of   lines of  ${\rm PG}_1(F v _i)$ by 
a copy of
the set of  blocks of  $\Delta _1$, 
subject only to the condition
\begin{itemize}
\item [$(\#)$]  there are distinct blocks,  neither of which is a line of {\bf  P},
whose span in  {\bf  P} is  ${\rm PG}_1(F v _i)$,
\end{itemize}
\noindent
and then apply all $g\in G$ to these sets of blocks in order to obtain the blocks in 
${\rm PG}_1((F v _i)^g),$ $g\in G; $ and
\item [(II)]   replace the set of     lines of   ${\rm PG}_1(F (v _i +\theta v_j))  $  by 
a copy of
the set of  blocks of 
$\Delta _2,$  
subject only to  $(\#)$,
 and then apply  all $g\in G$  to these sets of blocks in order to obtain the 
blocks in ${\rm PG}_1(F (v _i +\theta v_j) ^g),$ $g\in G. $
\end{enumerate} 
\indent
We need to check that these requirements can be met.

(i)  \emph{Satisfying} $(\#)$:
Let $\bar\Delta _s$ be an isomorphic copy of $\Delta _s$, $s=1$ or 2, whose set of points is
that of ${\rm PG}_1(F v)={\rm PG}_1(F v _i)$ or ${\rm PG}_1(F(v _i +\theta v_j))$.
Let $B_1$ and $B_2$ be any distinct blocks of $\bar\Delta _s$.  
Choose any permutation $\pi$ of the points of  ${\rm PG}_1(F v)$ such that the sets
$B_1^\pi$ and $B_2^\pi$ are not lines of ${\rm PG}_1(F v)$ and together
span ${\rm PG}_1(F v)$.
Using  $\bar\Delta _s^\pi$ in place of $ \bar\Delta _s$ satisfies $(\#)$.   (If $q+1\ge4$ then $B_2$ is not needed.)
 
(ii) 
\emph{These replacements are well-defined}:
For  (II), if $F (v _i +\theta v_j) ^g \cap F (v _i +\theta v_j) ^{g'}\ne 0 $ for  some $g,g'\in G,$
 then  $v _{i ^{g'}} +\theta v_{j ^{g'}}\in F (v _{i ^{g}} +\theta v_{j ^{g}}) $.
 Then either
 $v _{i ^{g'}} =v _{i ^{g}}$ and $v _{j ^{g'}} =v _{j ^{g}}$,    
 or $v _{i ^{g'}} =\alpha \theta v _{j ^{g}}$ and $\theta v _{j ^{g'}} =\alpha v _{i ^{g}}$
 for some $\alpha\in F^*$;
 but in the latter case we obtain  $1=  \alpha \theta $ and $ \theta=\alpha  $,
whereas  $\theta$ generates $F^*$. 
Thus,  $v _{i ^{g'}} =v _{i ^{g}}$, so the semiregularly of $G$ on $ \{ 1,\ldots , n\} $ implies that
 $g'=g$, as required.

It is trivial to see that ${\bf D} $ {\em is a design having the same parameters as
${\rm PG}_1(V)$.}  
Clearly   $G$ acts on the collection of subsets of $\mathfrak P$ occurring in (I) or (II):
we can view $G$ as a subgroup of both ${\rm Aut} {\bf D}$ and ${\rm PGL}(V)$.

We emphasize that the   sets in {\rm(I)}  and  {\rm(II)} occupy a tiny portion of the 
underlying projective  space:
most sets $Fv$  are unchanged.  More precisely,  in view of  the definition of
{\bf D}:
\begin{equation}
\label{block is line}   
 \begin{array}{llll}
 \mbox{\em Every block of {\,\bf  D} not contained  in a set {\rm(I)}  or  {\rm(II)} 
 is a line of {\,\bf  P}.}
 \\
  \mbox{\em Every line of {\,\bf  P} not  contained  in   set {\rm(I)}  or  {\rm(II)} 
 is a block of {\,\bf  D}.} 
 \end{array} 
 \end{equation}
 
Nevertheless, we will  distinguish between the   \emph{lines of  {\bf  P}}
and the \emph{blocks of {\bf D}}, even when the blocks happen to be lines.
A {\em subspace}  of {\bf D}  is a set of points  
that contains the block joining any pair of its points.  
(Examples:   (I) and (II) involve subspaces of {\bf D}.)
 A  {\em hyperplane} of {\bf D}
  is a  subspace  of {\bf D} that meets every block but does not contain every point. 
  We need further notation:
  \begin{equation}
   \label{block&line} 
 \mbox{
  \em Distinct $y,z\in \mathfrak P$ determine a block $yz$ of  {\,\bf  D} and a line 
   $\langle y,z\rangle$ of {\,\bf  P}. \hspace{8pt} }  \vspace{-4pt}
\end{equation}
  \begin{equation}
   \label{xyz} 
 \begin{array}{llll}
\mbox{\em For distinct $y,z\in \mathfrak P$ and $x\in   \mathfrak P- yz,$  }
\vspace{1pt}
    \\
\mbox{\em$\langle x|y,z\rangle =\bigcup\big\{xp\mid p\in y'z',  
 y'\in xy\hspace{-.5pt} -\hspace{-.5pt}\{x\},    z'\in xz\hspace{-.5pt} -\hspace{-.5pt}\{x\}, \hspace{-.5pt} \{y,z\}\hspace{-.5pt}\ne\hspace{-.5pt}  \{y',z'\}\big\}$.}  
      \end{array} 
\end{equation}
Here  (\ref{xyz}) depends only on {\bf  D} not on {\bf  P}, which will allow us to 
recover {\bf  P} from {\bf  D}.

       \begin{Lemma}
       \label{planes}
 If  $y,z\in \mathfrak P$ are distinct$,$   
 then there are more than $\frac{1}{2}|\mathfrak P|$   points $x\in  \mathfrak P-yz $  such that
   \begin{itemize}
  \item  [\rm(1)]   $\langle x,y,z\rangle$   is a plane
      of   {\,\bf  P}   every line  of  which$,$ except possibly $ \langle y,z\rangle,$  is a block of 
      {\,\bf  D}$,$ 
                 
 \item  [\rm(2)]   $ \langle x|y,z\rangle=\langle x,y,z\rangle ,$
 \item  [\rm(3)]  if $yz \subseteq  \langle x|y,z\rangle $ then  $  \langle y,z\rangle=yz ,$
    and                   
 \item  [\rm(4)] if $yz \not \subseteq  \langle x|y,z\rangle $ then   
 $\langle y,z\rangle$ is the union of the   pairs     
 $\{y_1,z_1\}\subset \langle x|y,z\rangle$  such that   $y_1z_1\not\subseteq  \langle x|y ,z \rangle.$
 \end{itemize}
   \end{Lemma}

\proof  
 Let   
\begin{equation}
\label{x}   
x\notin  yz\cup \mbox{$  \bigcup\hspace{1pt}$}\big \{ \langle  y,z,Fv\rangle  \mid 
\mbox{$Fv$ in (I) or (II)$ \big \}.$  }
\end{equation}
There are more  than
${(q^{4n}-1)/(q-1)} - n^2(q^6-1)/(q-1)  -(q+1)> \frac{1}{2}|\mathfrak P| $ such points $x$.
Clearly   $\langle x,y,z\rangle$  is a plane of {\bf  P}.

\smallskip
 (1)   Let $L\ne \langle y,z\rangle$ be a line of   $\langle x,y,z\rangle$,
 so   $\langle x,y,z\rangle = \langle y,z, L \rangle$.
 If $L$ is not a block of  {\bf  D}  then,  by  (\ref{block is line}),  $L$ is contained in some set $Fv$ in (I) or (II), so   $x\in \langle y,z,L\rangle\subseteq  \langle  y,z, Fv\rangle$ contradicts   (\ref{x}).

\smallskip
(2)    By (1), $\langle x,y\rangle $ and $ \langle x,z\rangle $ are blocks of {\bf  D}.
Let  $\{y',z'\}$ be as in  (\ref{xyz}).   Then  $\{y',z'\}\subset  \langle x,y,z\rangle$  
and $\langle y' ,z'\rangle \ne \langle y  ,z \rangle $.
By (1), $ y' z'=\langle y' ,z'\rangle \subseteq \langle x,y,z\rangle$ and
$xp = \langle x,p\rangle  \subseteq \langle x,y,z\rangle$ for each point $p$ of $\langle y' ,z'\rangle $.
Then $\langle x|y,z\rangle\subseteq \langle x,y,z\rangle$.
Each point of $\langle x,y,z\rangle$ lies in such a line $\langle x,p\rangle $;
since  that line is a block by  (1),    $\langle x,y,z\rangle\subseteq \langle x|y,z\rangle $. 
  
  \smallskip
(3)  If  $  yz\ne \langle y,z\rangle$  then,  by (\ref{block is line}),
 $  yz$  lies in some set $Fv$ in (I) or (II).
  By hypothesis   and (2),   $yz \subseteq \langle x|y,z\rangle\cap Fv=\langle x,y,z\rangle\cap Fv=\langle y,z\rangle$.  Thus,   $  yz= \langle y,z\rangle.$   
 
 \smallskip
 (4)   We have $yz\ne \langle y,z\rangle$  since 
 $ \langle y,z\rangle\subseteq  \langle x,y,z\rangle= \langle x|y,z\rangle$ by (2).
 By (\ref{block is line}), since $ \langle y,z\rangle$ is not a block it is contained in 
 some set $Fv$ in (I) or (II).
 
  For any  $\{y_1,z_1\}$
 in (4) we have           $\{y_1,z_1\}\subseteq \langle x|y,z\rangle= \langle x,y,z\rangle$    by (2),
and $y_1z_1\not\subseteq  \langle x,y,z\rangle$,
so  $\langle y_1,z_1\rangle$ is not a  block of {\bf  D} and hence 
$\langle y_1,z_1\rangle= \langle y,z\rangle$ by (1).

On the other hand,  consider an arbitrary pair
    $\{y_1,z_1\} \subset \langle y,z\rangle \subset Fv$.  Then 
$y_1 z_1 \subset Fv$ by the definition of {\bf  D}.  Since  $\langle y,z\rangle $ is not a block,
$y_1 z_1\not  \subseteq \langle y,z\rangle=\langle x|y,z\rangle\cap Fv$  by (2),
so $y_1 z_1\not  \subseteq \langle x|y,z\rangle$.
Thus,  $\langle y,z\rangle$ is the union of the  pairs     $\{y_1,z_1\}$ in (4). \qed    
    
      \smallskip  \smallskip 
           \noindent
\emph{Proof of} Theorem~\ref{design construction3}(i).
  We first recover the lines of {\bf  P} from {\bf  D}.  
  For distinct  $y,z\in \mathfrak P$,
use each $x\notin yz$  in Lemma~\ref{planes}(3) or (4)   
in order to obtain,  more than $ \frac{1}{2}|\mathfrak P| $ times,
the same   set  of points that must be $\langle  y,z\rangle $.
    
\emph{We have now   reconstructed all lines of  {\,\bf  P} as subsets of $\mathfrak P$}.
 Then we  have    also recovered   ${\bf  P},$ $V,$  ${\rm \upGamma L}(V)$ 
  and   ${\rm P\upGamma L}(V)$,   so that
${\rm Aut} {\bf D}$ is induced by a subgroup of~${\rm Aut} {\bf P}={\rm P\upGamma L}(V)$. 

Any block of {\bf D}  that is not a line of {\bf P}
spans a   2-space or 3-space of {\bf P} occurring in
 some 3-space  ${\rm PG}_1(F v)$ in (I) or (II),
and spans at least a   4-space of {\bf P} together with any  block  in any 
${\rm PG}_1(F v')\ne  {\rm PG}_1(F v)$.  
Any two  blocks of {\bf D}   that    are not  lines of {\bf P} and lie in the same set  in (I) or (II)
span at most a 3-space of {\bf P};  by   $(\#)$  each  set  in (I) or (II) is spanned by two such blocks.
  
\emph{This  recovers all subsets {\rm(I)}  and  {\rm(II)} of $\,\mathfrak P$}
from {\bf D} and {\bf P}.  Moreover, 
the fact that $\Delta_1\not\cong\Delta_2$ specifies which of these  subspaces of {\bf D}
have  type (I) (or (II)).

We next determine the $F$-structure of $V$ using {\bf D}.
We claim that \emph{the subgroup of ${ {\rm PGL}}(V)$ fixing each  set in 
{\rm(I)} or {\rm(II)} arises from 
 scalar multiplications by members of $F^*$.}  
Clearly  such scalar multiplications behave  this way. Let $h\in { {\rm PGL}}(V)$
behave as stated.  Let $\hat h\in { {\rm GL}}(V)$ induce $h$.
Then
 $\hat h\colon xv_i\mapsto (xA_i)v_i$ for  all $x\in F$  and a $4\times4$  
invertible matrix $A_i$ over $K$.
If $ij$ is an edge of $\Gamma$ and $x\in F$, then 
$(x(v_i+\theta v_j))^{\hat h}= (xA_i )  v_i+ ((x\theta)A_j ) v_j$ is in $F(v_i+\theta v_j)$, so
$(xA_i ) \theta = (x\theta)A_j  $. Since $ij=ji$,  also $(xA_j ) \theta = (x\theta)A_i $,  so
$(x\theta\theta)A_i=((x\theta)A_j)\theta=   (xA_i) \theta\theta$, and $A_i$ commutes with multiplication
by $\theta^2$.  By Schur's Lemma, $xA_i =   xa_i  $ 
for  all $x\in F$  and some 
$a_i  \in F^*$.
Then $xa_i\theta=x\theta a_j $, 
so   $a_i=a_j$.  
Since  $\Gamma $ is connected, 
all $a_i$ are equal,  proving our claim.

In particular,      the field  $F$ and   the $F$-space  $V_F$ 
can be reconstructed from    {\bf D}.  
Since {\rm Aut}{\bf D}  normalizes $F^*$ by the preceding paragraph,   
$ {\rm Aut} {\bf D}\le{\rm P\upGamma L}(V_F).$
We know that $G$  is inside both ${\rm Aut} {\bf D}$ and ${\rm PGL}(V)$.
Since the sets in (II) correspond to (ordered) edges of $\Gamma$,
$ {\rm Aut} {\bf D}$ induces  $ {\rm Aut} \Gamma\cong G$ on  the collection  of sets in (I).

Let $h\in {\rm Aut} {\bf D}\le {\rm P\upGamma L}(V_F)$.  Multiply $h$  by an element of $G$ in order to have~$h$ fix  
all $Fv_i$.  Let $\hat h\in { {\rm \upGamma L}}(V_F)$ induce $h$,
with associated field automorphism $\sigma\in {\rm Aut} F $. 
For each $i$ we have  $v_i^{\hat h}=a_iv_i$ for some ${a_i\in F^*}$. 
Let $ij$ be an edge of  $\Gamma$ and write $b=a_j/a_i$.   As above,   
${F (v _i +\theta v_j)^{\hat h}}= {F (a_i v _i +\theta ^\sigma a_jv_j)}= 
{F(v _i +  \theta ^{\sigma} b v_j) }$~and $F (\theta v _i +v_j)^{\hat h}
= {F (\theta ^\sigma a_iv _i +a_jv_j)}
= F ( v _i +\theta ^{-\sigma }b v_j)$ 
both have type (II), so  
$ \theta^\sigma b=  \theta^{\pm1}  $  and $ \theta^{-\sigma}b=  \theta^{\mp1} .$
Then ${b^2=1}$, $  \theta^\sigma =  \pm\theta^{\pm1}  $, and hence $\sigma=1$ and $b=1$
since $\theta$ generates $F^*$. 
The  connectedness of $\Gamma $   implies that
 all $a_i$  are equal:  $\hat h$ is scalar multiplication by $a_1$.
 
Since $h$ fixes $Fv_1$ it induces an automorphism of the subspace of  {\bf D}
determined by $Fv_1$.
By (I) and our condition  
on  $\Delta _1$, $h$ fixes a point $K cv_1  $  of $Fv_1$.  
Then $K cv_1  =  (Kcv_1) ^h= K ca_1v_1 $, so
$a_1\in K$.  Thus,
$h=1$  on ${\mathfrak P}$  and  $ {\rm Aut} {\bf D} \cong G$.\qed 

\section{A simpler projective construction}
\label{simpler projective}
We need a fairly weak result (Proposition~\ref{needed designs}) concerning designs with the parameters of  ${\rm PG}_ 1(3,q)$.  
We know of two published constructions for designs having those parameters, 
due to Skolem   \cite[p.~268]{Wi}
and  Lorimer \cite{Lo}.  However, isomorphism questions seem difficult using
their descriptions.  Instead, we will use a method that
imitates   \cite{Sh,Ka} (but which was hinted at by Skolem's idea).  
  
Consider a hyperplane $X$ of  ${\bf P}={\rm PG}(d,q)$, $d\ge3$;
we identify {\bf P} with  ${\rm PG}_1(d,q)$.
Let $\pi $ be any permutation of the  points of $X$. Define a geometry ${\bf D}_\pi $ as follows:
\vspace{-2pt}
\begin{itemize}
\item [] the set $\mathfrak P$  of points is the set of  points of {\bf P},  and
\item [] 
 blocks are of two sorts:
 \begin{itemize}
 \item [] the lines of {\bf P}  not in $X$, and
\item []  the sets $L^\pi $ for lines $L \subset X$.
\end{itemize}
\end{itemize}
\vspace{-2pt}

Once again it is trivial to see that ${\bf D}_\pi $   
is a design having the same parameters as {\bf P}.  
Note~that $\pi $ has nothing to do with the incidences
between points and the blocks not in~$X$.  

We have  a  hyperplane $X$ of  ${\bf D}_\pi $ such that the  blocks of ${\bf D}_\pi $
 not in $X$ are   lines of a projective space {\bf P}  for which  $\mathfrak P$   is the set of points.
 We claim that  \emph{the lines of  this projective space can be recovered 
  from ${\bf D}_\pi $ and $X$.}
   Namely,  we have all points and lines of   {\bf P}
   not in $X$.  For distinct $y,z\in X$ and $x\notin X$, 
 the set   $\langle x|y,z\rangle$  in  (\ref{xyz}) consists of the points of  the plane  $\langle  x,y,z  \rangle$
 of   {\bf P}, and  $\langle x|y,z\rangle \cap X$ is the line $\langle y,z\rangle$.
   We have now obtained all lines of the original projective  space~{\bf P}, as claimed. 
   It follows that ${\rm Aut}{\bf D}_\pi \le {\rm Aut}{\bf P}.$

The symbol $X$ is ambiguous:  it will now mean either a set of points or
a hyperplane of the underlying \emph{projective space}  (as in the next result).  It will not refer to $X$ together with a different set of lines produced by a permutation $\pi$.    

\begin{prop}
\label{double cosets} 
The designs
${\bf D}_\pi $ and ${\bf D}_{\pi '}$ are isomorphic by an isomorphism sending $X$ to itself if and only if ${\pi }$ and ${\pi '}$ are in the same 
${\rm P\upGamma L}(X), {\rm P\upGamma L}(X)$ double coset in ${\rm Sym}(X)$.

Moreover$,$  the pointwise stabilizer of $X$ in  ${\rm Aut}{\bf D}_\pi $ is transitive on the points outside of $X,$   and the stabilizer $({\rm Aut}{\bf D}_\pi )_X$  of $X$ 
induces ${\rm P\upGamma L}(X)\cap {\rm P\upGamma L}(X)^\pi $ on~$X$.%
\end{prop}
\proof
Let $g{\colon\!} {\bf D}_\pi \to{\bf D}_{\pi '}$ be such an isomorphism.  We just saw that 
{\bf P} is naturally reconstructible from either design.  It follows that 
$g$ is   a collineation of {\bf P};  its restriction $\bar g$ to $X$ is in  ${\rm P\upGamma L}(X)$.

If $L\subset  X$ is a  line of {\bf P} then $g$   sends the block 
$L^\pi \subset  X$ of  ${\bf D}_\pi $   to a block   $L^\pi {}^g\subset  X$  of  ${\bf D}_{\pi '}$.
Then  $L^\pi {}^g{}^{\pi '^{-1}}$ is a line of  {\bf P}, so that 
 $\pi \bar g \pi '^{-1}$ is a permutation of the points of the hyperplane $X$ of  {\bf P}
 sending lines to lines, and hence is an element $  h\in {\rm P\upGamma L}(X)$.
 Thus,  $\pi $ and $\pi '$ are in the same ${\rm P\upGamma L}(X), {\rm P\upGamma L}(X)$ double coset.
 
 Conversely, if $\pi $ and $\pi '$ are  in the same ${\rm P\upGamma L}(X), {\rm P\upGamma L}(X)$ double coset let  $\bar g,   h\in  {\rm P\upGamma L}(X)$ with $\pi \bar g\pi '^{-1}=  h $.
 Extend  $\bar g$ to $g\in {\rm Aut} {\bf P} $ in any way.  We claim that $g$ is an isomorphism
  $ {\bf D}_\pi \to{\bf D}_{\pi '}$.  It  preserves incidences between
  blocks not in $X$ and points of {\bf P}
  since $g\in {\rm Aut} {\bf P} $ and those incidences have nothing to do with $\pi$ and $\pi'$.  
  Consider an incidence $x\in B\subset X$
  for a block $B$ of ${\bf D}_\pi $.  Then $B = L^\pi $ for a line $L\subset X$.
  Since $g\in {\rm Aut} {\bf P} $,  $x^g\in B^g=B^{\bar g}=L^{\pi \bar g}=(L^{  h})^{\pi '}$,
  which  is a block of ${\bf D}_{\pi '}$,   as required.
  
  For the final assertion, the  pointwise stabilizer  of  $X$ in ${\rm Aut}{\bf P}$
  is in  ${\rm Aut}{\bf D}_\pi $ by the definition of  $ {\bf D}_\pi $.
 We have seen that the group induced   on $ X$ by 
  ${\rm Aut}{\bf D}_\pi $   corresponds to the pairs  $(\bar g ,  h )\in
  {\rm P\upGamma L}(X)\times {\rm P\upGamma L}(X)$ satisfying  $\pi \bar g\pi ^{-1}=  h $.
   \qed
 \smallskip\smallskip
 
 Note that there are many extensions $g$ of $\bar g$ since the designs
 ${\bf D}_\pi $ have many automorphisms inducing the identity on $X$.
 Double cosets arise naturally in this  type  of result; compare \cite[Theorem 4.4]{Ka}.

Let  $v_i=(q^i -1)/(q-1)$.   

\begin{coro}
\label{numbers}
There are at least  $v_d!/ v_{d+1}|{\rm P\upGamma L}(d,q)|^2$
pairwise nonisomorphic designs having the same parameters as {\bf P}.
\end{coro}
\proof
Fix $\pi $ in the proposition.  
There are  at most $v_{d+1}$ hyperplanes $Y$ of ${\bf D}_\pi $
(as in   \cite[Theorem 2.2]{JT}).
By the proposition there are then at most $|{\rm P\upGamma L}(X)|^2$ choices for
  $\pi ' $ such that   ${\bf D}_\pi \cong{\bf D}_{\pi '}$ by an isomorphism sending $Y$ to $X$.  
Since there are $v_d!$  choices for  $\pi $ we obtain the stated lower bound.\qed

 \rema   \em
 We describe a useful trick.  \em
  A  transposition $\sigma$ and a $3$-cycle  $\tau$ are in different
${\rm P\upGamma L}(d,q), {\rm P\upGamma L}(d,q)$ double cosets in ${\rm Sym}(N),$ 
$ N=   (q^d-1)/(q-1),$ if  $d\ge 3$ and we exclude the case $d=3,q=2$.  \rm 
For, if $\sigma  g =  h \tau$ with $g,h\in {\rm P\upGamma L}(d,q)$  then 
$ g^{-1}  h =g^{-1}\cdot \sigma g \tau ^{-1}   =
{\sigma} ^{  g}{\tau} ^{-1} \in {\rm P\upGamma L}(d,q)$ fixes at least $N-5$ points,
and hence is 1 by our restriction on  $d$, whereas ${\sigma} ^{  g}\ne {\tau} $.

\begin{prop}
\label{needed designs}
For any $q$ there  are  two  designs having the parameters of ${\bf P}={\rm PG}_1(3,q)$ and not isomorphic to one another  or  to
${\bf P} ,$ for   one of which the automorphism group fixes a point.
\end{prop}
\proof
If $q=2 $ then there are even such designs   with  trivial automorphism group \cite{CCW}.
(Undoubtedly     such designs  exist for all $q$.)

Assume that $q>2$.  The    preceding corollary and remark provide us with   two nonisomorphic designs.   It remains to deal with the final assertion constructively. 

Let $\pi $ be a transposition $(x_1,x_2)$ of $X$.  We will show that  ${\bf D}_\pi $ behaves as stated.

First note that each $g\in  {\rm Aut}{\bf D}_\pi \le {\rm Aut}{\bf P}$ fixes $X$.
For, suppose that $Y=X^g\ne X$ for some $g$.  The blocks in $Y$ not in $X$ are lines
of {\bf P}.  Then the same is true of the 
 blocks in $Y^{g^{-1}}=X$ not in $X^{g^{-1}}$.  This contradicts the fact that $\pi $ sends all lines $\ne\langle x_1,x_2\rangle$  of  {\bf P} inside $X$ and on $x$  to sets that are not lines of {\bf P}.
      
By Proposition~\ref{double cosets},  
${\rm Aut}{\bf D}_\pi =({\rm Aut}{\bf D}_\pi)_X  $  
 induces  ${\rm P\upGamma L}(X)\cap {\rm P\upGamma L}(X)^\pi $ on $X$.  
Let $\pi \bar g\pi ^{-1}=  h $ for $\bar g,  h\in {\rm P\upGamma L}(X)$.
Then $\bar g^{-1}  h  =\pi ^{\bar g} \pi ^{-1}$
 is a collineation of $X$ that moves at most $2\cdot 2$ points of $X$
 and hence fixes at least $(q^2+q+1) -2\cdot 2 >  q+\sqrt q +1 $ points. 
 By elementary (semi)linear algebra,     the only such collineation is 1,
so that  $\bar g=  h$ commutes with $\pi $ and hence   fixes the line $\langle x_1,x_2\rangle$.
Then $\bar g$ also fixes a point of $X$ and hence of  ${\bf D}_\pi $.  
\qed
            
\rema\rm  By excluding   the possibilities $q\le8$ and $q$ prime  in the previous section 
we could have used nondesarguesian projective planes  (and  ${[F\colon\!K]=3}$).

   \section{A simple affine construction}
  \label{affine construction}
  \emph{We now consider }Theorem~\ref{design construction3}(ii).
    The proof is similar to that of   Theorem~\ref{design construction3}(i).  
That result handles the cases $q=3,$ 4 or 5,   but we  ignore this  and only assume that $q>2$.
  
Let $G$ and $\Gamma$ be as in Section~\ref{Construction}. 
This time we use  $K={\bf F}_q\subset F={\bf F}_{q^3}$; once again $\theta$ generates $F^*$.  
Let $V_F$ be an $n$-dimensional  vector space over   $F$,
with basis  $  v_1,\ldots ,v_n $.
View $V_F$  as a vector space $V$ over $K$.  
  If $Y$ is a set of points of {\bf  A} then $\langle Y\rangle$ denotes the smallest affine subspace
    containing $Y$.

We will modify the 
point-line design $ {\rm AG}_1(V)$ of  ${\bf A}= {\rm AG}(V)$, using  nonisomorphic designs 
$  \Delta _1,\Delta _2 $ 
whose parameters are those of   $ {\rm AG}_1(3,q)$  but are not isomorphic to 
that design, chosen so that   ${\rm Aut} \Delta _1$ fixes at least two   points
(Proposition~\ref{needed resolvable designs}). 

Our design $\bf D$ has $V$ as its set of points.  Most   blocks of $\bf D$  
are lines of  {\bf A}, with   exceptions involving the
sets $Fv, $ $0\ne v\in V$,   in Section~\ref{Construction}(I,\,II),
where now {\em $Fv$ is viewed as a $3$-dimensional affine space.}   

As before,  the set of lines of  
 $ {\rm AG}_1(Fv_i)$  or  $ {\rm AG}_1(F (v _i +\theta v_j))$ 
 is replaced by a copy of the set of blocks of $\Delta_1$ or   $\Delta_2$.
 This time, for each of these  we require 
 \begin{itemize}
\item [$(\#')$]  there are distinct blocks,  each of which spans a plane of  {\bf  A},
such that  the intersection of those planes is a line.
\end{itemize} 
Clearly, these two blocks span a 3-space.
(When $q>3$ it would be marginally easier to require that there is a single block 
that spans a 3-space.) 
 Condition $(\#')$ can be satisfied
  exactly as in   \emph{Satisfying} $(\#)$ in Section~\ref{Construction}.
Since different sets $Fv  $ meet  only in a single point, the modifications made 
 inside  them  are unrelated.
Once again it is easy to check that this produces a design  {\bf D} 
with the desired parameters for which $G\le {\rm Aut} {\bf D}$. 

As in Section~\ref{Construction},   
most sets $Fv$  are unchanged.
In view of  the definition of  {\bf D}, the analogue of (\ref{block is line}) holds.
We use the natural analogues of   definitions  (\ref{block&line})  and  (\ref{xyz}),
using {\bf A} in place of {\bf P} and $V$ in place of $\mathfrak P$.

       \begin{Lemma}
       \label{aplanes}
 If  $y,z\in V$ are distinct$,$   
 then there are more than $\frac{1}{2}|V|$   points $x\in V-yz $  such that
   \begin{itemize}
  \item  [\rm(1)]   every line  of the plane
      $\langle x,y,z\rangle$   of  {\bf A}$,$ except possibly $ \langle y,z\rangle,$  is a block of {\bf  D}$,$ 
                 
 \item  [\rm(2)]   $ \langle x|y,z\rangle =\langle x,y,z\rangle ,$
 \item  [\rm(3)]  if $yz \subseteq  \langle x|y,z\rangle $ then  $  \langle y,z\rangle=yz ,$
    and                   
 \item  [\rm(4)] if $yz \not \subseteq  \langle x|y,z\rangle $ then   
 $\langle y,z\rangle$ is the union of the   pairs     
 $\{y_1,z_1\}\subset \langle x|y,z\rangle $  such that   $y_1z_1\not\subseteq  \langle x|y ,z \rangle  .$
 \end{itemize}
 \end{Lemma}

  \proof
  Using $x$ in (\ref{x}),  this is proved exactly as in Lemma~\ref{planes} except for (2), where 
  we need to consider parallel  lines using blocks that are lines by (1).  
  Clearly $ \langle x|y,z\rangle \subseteq \langle x,y,z  \rangle$;  
  we must show that  $ \langle x,y,z\rangle \subseteq \langle x|y,z\rangle.$
  In (\ref{xyz}),  for $p$ in the line 
  $y'  z ' =\langle y' ,z '\rangle$ of $\langle x,y,z\rangle$  parallel to  $\langle y,z\rangle$,
   the blocks $xp\subset \langle x|y,z\rangle$ cover all points 
   of  the plane $\langle x,y,z\rangle$  except for those in the line $L$ on $x$ parallel 
  to  $\langle y,z\rangle$.   If $y'\in xy-\{x,y\}$ and 
  $p'=y'z\cap L$, then $L=xp'\subset \langle x|y,z\rangle $,
  so $ \langle x,y,z\rangle \subseteq \langle x|y,z\rangle.$  \qed
    
\smallskip  \smallskip    
\noindent\emph{Proof of\,\,\,}Theorem~\ref{design construction3}(ii). 
First  recover all lines of  {\bf  A}  from {\bf D} exactly 
as in the proof of  Theorem~\ref{design construction3}(i).
This also produces  both  the $K$-space $V$ and ${\rm A\upGamma L}(V)$ from {\bf D}.

We   recover all  subsets (I) and  (II) essentially as before.     
 Consider a pair $B,B'$ of blocks of  {\bf D}  behaving as in $(\#')$:  $\langle B\rangle $
 and  $\langle B'\rangle $ are planes and $\langle B\rangle \cap \langle B'\rangle $ is a line.
 Since distinct subsets in  (I) or  (II) do not have a common line, each such pair 
 $B,B'$ spans a  subset  in  (I) or  (II).
   Thus, by $(\#')$ we have obtained  each subset  in  (I) or  (II)   from  {\bf D} and {\bf A} using
some pair $B,B'$.   Once again, the fact that $\Delta_1\not\cong\Delta_2$ 
specifies which of these  subspaces of {\bf D}    have  type (I) (or (II)).

The subsets (I) all contain $0$,   and  ${\rm Aut} {\bf D}$ fixes their intersection,
so     ${\rm Aut} {\bf D}$ is induced by a subgroup of ${\rm A\upGamma L}(V)_0
 ={ {\rm \upGamma L}}(V )$.  
  
 Recover the field $F$ exactly as in the proof of  Theorem~\ref{design construction3}(i).
Once again,  $ {\rm Aut} {\bf D}$ is a subgroup of ${\rm \upGamma L}(V_F)$ that
induces  $ {\rm Aut} \Gamma\cong G$ on  the collection  of sets  in (I).   

By repeating the argument  at the end of the proof of  Theorem~\ref{design construction3}(i)  we 
reduce to the case of  $ h\in {\rm Aut} {\bf D}$ 
fixing all sets  in (I) and  acting  on $V$  as $v\mapsto av$ for some 
$a\in F^*$.   We chose $\Delta_1$ so that $ {\rm Aut} \Delta_1$   fixes at least two  of its points.   It follows that  $a=1$, so that $h=1$ and  $ {\rm Aut} {\bf D}\cong G$.   \qed
  
        \section{A simpler affine construction}  
      \label{simpler affine}  
  Consider a plane $X$ of 
${\bf A}= {\rm AG}(3,q)= {\rm AG}(V)$, $q>2$; we identify {\bf A} with  $ {\rm AG}_1(3,q)$.
Let $\pi $ be any permutation of the  points of $X$.
Define a geometry ${\bf D}_\pi $ as follows:
\vspace{-2pt}
\begin{itemize}
\item [] the set $V$  of points is the set of  points of {\bf A},  and
\item [] 
 blocks are of two sorts:
 \begin{itemize}
 \item [] the lines of {\bf A} not in $X$, and
\item []  the sets $L^\pi $ for lines $L\subset X$.
\end{itemize}
\end{itemize}
\vspace{-2pt}

Once again it is trivial to see that ${\bf D}_\pi $   is a design having the same parameters as {\bf A}.    

As in Section~\ref{simpler projective},  the  blocks of ${\bf D}_\pi $
 not in $X$ are   lines of an affine space {\bf A} for which  $V$   is the set of points.
As in Sections~\ref{simpler projective} and  \ref{affine construction}, \emph{the lines of  this affine space can be recovered 
 from ${\bf D}_\pi $}  using the analogue of   (\ref{xyz}).

\begin{prop}
\label{double cosets 2}
The designs
${\bf D}_\pi $ and ${\bf D}_{\pi '}$ are isomorphic by an isomorphism sending $X$ to itself if and only if ${\pi }$ and ${\pi '}$ are in the same 
${\rm A\upGamma L}(X), {\rm A\upGamma L}(X)$ double coset in ${\rm Sym}(X)$.
This produces at least $q^2!/q(q^2+q+1)|  {\rm A\upGamma L}(2,q)  | ^2$ pairwise 
nonisomorphic designs having the same parameters as   ${\rm AG}_1(3,q)   $.

Moreover$,$  the pointwise stabilizer of $X$ in  ${\rm Aut}{\bf D}_\pi $ is transitive on the points outside of $X,$   
and $({\rm Aut}{\bf D}_\pi)_X $  
induces ${\rm A\upGamma L}(X)\cap {\rm A\upGamma L}(X)^\pi $ on  $X$.
\end{prop}
\proof
This is the same as for Proposition~\ref{double cosets} and Corollary~\ref{numbers}.
   \qed  

\begin{prop}
\label{needed resolvable designs}
For any $q  \ge3$ there  are at least  two designs having the parameters of ${\bf A}= {\rm AG}_1(3,q),$
   not isomorphic to one another  or  to
${\bf A}  ,$    such that the automorphism group of one of them fixes at least two points.
\end{prop}
\proof
The bound in the preceding proposition provides us with many nonisomorphic designs.
We need to deal with the requirement concerning    automorphism groups.
By \cite{LR} we may assume that $q\ge4$.

Let $\pi\in {\rm Sym}(X)$   be a 4-cycle $ (x,x_1,x_2,x_3)$, where
$x_1,x_2,x_3$ are on a line not containing $x$.
We will show that ${\bf D}_{\pi}$ behaves as required.

Let $g\in  {\rm Aut}{\bf D}_{\pi  }$.  As  in the proof of Proposition~\ref{needed designs},
   $g$ fixes $X$  and induces a collineation $\bar g$ of the subspace $X$  of {\bf A}.
 By Proposition~\ref{double cosets 2},  $\pi  \bar g=  h \pi $ with $\bar g,  h\in {\rm A\upGamma L}(X)$.    
As before,    $\bar g^{-1}  h =\pi  ^{\bar g}\pi  ^{-1}$ is a collineation of $X$ that 
 fixes at least $q^2 -2\cdot 4>  q  $ points  as $q  \ge  4$.~Then  
  $\bar g=  h$   and $ \pi  ^{\bar g}=\pi  $.  
 Since    the collineation $\bar g$  commutes~with~$\pi $
 it fixes $ \{x,x_1,x_2,x_3\}$ and hence also $x$, and so
   is the identity on the support of $\pi $.  Thus, $ {\rm Aut}{\bf D}_{\pi  }$ is the identity on that support.  \qed%

        \section{Steiner quadruple systems}  
      \label{Steiner quadruple systems}  
We have avoided ${\rm AG}(d,2)$ in the preceding two sections.   
Here we briefly comment about those spaces in the context of $3$--$(v,4,1)$-designs 
 (Steiner quadruple systems), outlining a proof of the  following result in  \cite{Me}.
 
 \begin{theorem}
 \label{SQS}
If $G$ is a finite group  then there are   infinitely many integers $v $ such that there is a
$3$--$(v,4,1)$-design ${\bf D}$ for which  ${\rm Aut} {\bf D}\cong G$.
 \end{theorem}
 \proof
 Let $K={\bf F}_2\subset F= {\bf F}_{16}$ be as in Section~\ref{Construction},  with $\theta $
 a generator of $F^*$.  Let $V_F $ be a vector space with basis 
 $v_1,\dots,v_n$, viewed as a $K$-space $V$.  This time we modify the 3-design 
 ${\rm AG}_2(V)$ of points and (affine) planes of   $V$.  We use
 nonisomorphic designs $\Delta_1, \Delta_2$ having the parameters of ${\rm AG}_2(4,2)$ 
 but not isomorphic to that design, and such that  ${\rm Aut}(\Delta_1)=1$ \cite{KOP}.
 
Once again our design $\bf D$ has $V$ as its set of points.  Most   blocks of $\bf D$  
are planes of  {\bf A}, with   exceptions involving the
sets $Fv, $ $0\ne v\in V$,   in Section~\ref{Construction}(I,\,II),
where now {\em $Fv$ is viewed as a $4$-dimensional affine space.}   
As before,  the set of planes of 
 $ {\rm AG}_2(Fv_i)$  or  $ {\rm AG}_2(F (v _i +\theta v_j))$ 
 is replaced by a copy of the set of blocks of $\Delta_1$ or   $\Delta_2$.
 This time, for each of these  we require 
   \begin{itemize}
\item [$(\#'')$]  there are distinct blocks,  each of which spans a 3-space of  {\bf  A},
such that  the intersection of those  3-spaces is a plane.
\end{itemize}  
 Once again it is easy to check that this produces a design  {\bf D} 
with the desired parameters for which $G\le {\rm Aut} {\bf D}$. 
  
 {Distinct $x,y,z\in V$ determine a block $xyz$ of  {\,\bf  D} and a plane 
   $\langle x,y,z\rangle$ of {\,\bf  A}. }
   For distinct $x,y,z$ and $w\notin xyz,$  instead of (\ref{xyz}) we  use
 {$\langle w|x,y,z\rangle\!=\bigcup\big\{abc\mid
 a\in wxy-\{w\}, b\in wxz-\{w\}, c\in wyz-\{w\}$, with $a,b,c$ distinct and not all in   $\{x,y,z\}\big\}$.}
    
  As before, all planes of {\bf A} can be recovered from {\bf D}, this time using 
  various sets $\langle w|x,y,z\rangle$.   Also the sets in (I) and (II) can be recovered, 
  as can   $F$, and    
  the argument at the end of Section~\ref{affine construction} goes through as before. \qed
  
        \section{Concluding remarks}  
      \label{Concluding remarks}  
         \Remark\rm
          \label{additional}  
 When considering possible consequences of this paper it became clear that additional properties of  our designs    should   also  be mentioned.   
   \smallskip  \smallskip

 \begin{enumerate}
\item [\rm(1)]
Additional properties of the design {\bf D}  
in Theorem~\ref{design construction3}(i).

\begin{enumerate}
\item[\rm(a)]
{\em ${\rm PG}(3,q)$-connectedness}.
The following graph is  connected: the vertices are the subspaces of {\bf D} 
isomorphic to ${\rm PG}_1(3,q)$, with two joined when they meet.
\item[\rm(b)]
{\em ${\rm PG}(n-1,q)$ generation}.
{\bf D} is generated by its subspaces isomorphic to ${\rm PG}_1(n-1 ,q)$.

\item[\rm(c)]
Every point of {\bf D} is in a  subspace isomorphic to
${\rm PG}_1(n-1 ,q)$ (in fact, many of these).

\item[\rm(d)]  More than $q^n$ points are moved by every nontrivial automorphism of {\bf D}.

\end{enumerate}
\item [\rm(2)]
Additional properties of  the design {\bf D}  in Theorem~\ref{design construction3}(ii).

\begin{enumerate} 
\item[\rm(a)]
{\em ${\rm AG}(3,q)$-connectedness}.
The following graph is  connected: the vertices are the subspaces of {\bf D} 
isomorphic to ${\rm AG}_1(3,q)$, with two joined when they meet.
\item[\rm(b)]
{\em ${\rm AG}(n ,q)$ generation}.
{\bf D} is generated by its  subspaces isomorphic to ${\rm AG}_1(n ,q)$.

\item[\rm(c)]
Every point of {\bf D} is in a  subspace isomorphic to
${\rm AG}_1(n ,q)$   (in fact, many of these).

\item[\rm(d)]  More than $q^n$ points are moved by every nontrivial automorphism of {\bf D}.

 \end{enumerate}
 
\item [\rm(3)]
Additional properties of  the design {\bf D}  
in Theorem~\ref{SQS}.  This time versions of (2a) (using ${\rm AG}_2(4,2)$-connectedness),
(2b), (2c)  and (2d) hold.
\end{enumerate}
   \smallskip  \smallskip
   
These reflect the fact that the sets of points in (I) or (II) cover a tiny portion of the 
underlying projective or affine space:   a subset of the points determined by $F$-linear 
combinations of at most two of the $v_i$.   
For (1a),  it is easy to see that any point in  $\mathfrak P$ lies in a 4-space of $V$
that contains some point $ K \beta \sum_i\!v_i$, $\beta\in F^*$,
 and meets each set  in (I) or (II) in at most a point; by (\ref{block is line})  this produces 
 a subspace of {\bf D} isomorphic to ${\rm PG}_1(3,q)$. Moreover, 
   all $K \beta \sum_iv_i$ lie in   $ F( \sum_i\!v_i)$,  which also produces 
    a subspace of {\bf D} isomorphic to ${\rm PG}_1(3,q)$. 
 
For (1b) we give  examples of subspaces of  $V$:

\noindent
$\langle v_1 +\theta ^2v_2,  v_2+\theta ^2v_3+\theta ^iv_4, \dots ,v_{n-2}+\theta ^2v_{n-1}+\theta^i v_n,  v_1+v_2+v_4+v_5, \theta(v_1+v_2+v_4+v_5)\rangle$

\noindent
 for $2< i<q^4-1$.  Each of these misses all sets in (I) or (II), and hence 
  determines  a  subspace of {\bf D} isomorphic to
${\rm PG}_1(n-1 ,q)$. These subspaces  generate a subspace of {\bf D} containing
the points  $K(\theta^i-\theta^3)v_n$, $3< i<q^4-1$,  and hence also 
 ${\rm PG}_1(Fv_n)$.   Now  permute the subscripts to generate {\bf D}. 
 
Part (1c) holds by using $K$-subspaces similar to the above ones. 
 There are clearly projective spaces of larger dimension that are subdesigns of {\bf D}.
 
 Part (1d) depends on the semiregularity of $G$ on $\{ v_1,\ldots ,v_n \}$.  Use the points
 $K  \sum_i\alpha _i v_i$ with $\alpha_1=1$ and  $\alpha_i\in F-\{1\}$ for $i>1$, 
 where each $\alpha \in F-\{1\}$ occurs either for 0 or at least  two basis vectors  $v_i$. 
 The lower bound $q^n$ is   easy to obtain  but very poor.

Both (2) and (3) are handled as in (1).

\Remark \rm
In (II) we used the $K$-subspaces  $F (v _i +\theta v_j)$.  We could have used   
subspaces  $F (v _i +\theta_rv_j)$, $r=1,\dots,s,$ for various  $\theta_r$, together with
 further nonisomorphic designs $\Delta_{2,r}$
(which are needed  to distinguish among the $F (v _i +\theta_r v_j)$).   All proofs go through without difficulty, 
as do the additional properties in the preceding remark.

\Remark \rm
Each of our designs    has the same parameters as 
some  $ {\rm PG}_1(V)$  or  $ {\rm AG}_1(V)$.    
What is needed is a  much  better type of result, such as:   
\emph{for each finite group $G$ there is an integer $f(|G|)$ such that$,$ 
if $q$ is a prime power and if $v>f(|G|)$ satisfies the necessary conditions for
the existence of a $2$--$(v,q+1,1)$-design$,$ then there is such a design  {\bf D}
for which ${\rm Aut}{\bf D}\cong G$.}
When $q=2$ this result is proved  in a  sequel to the present paper \cite{DoK}.

          \smallskip  \smallskip   \smallskip  \smallskip
{\noindent\em Acknowledgements.} \rm
I am   grateful to Jean Doyen for providing me with a clear description of  Skolem's construction and for helpful  comments concerning this research.
I am also  grateful to  a referee for many helpful comments. 
This research was supported in part by a grant from the Simons Foundation.

\end{document}